\pdfoutput=1
\documentclass[12pt]{article}
\RequirePackage{filecontents}

\usepackage{geometry}
\usepackage{youngtab}
\usepackage{tikz}
\usepackage{placeins}
\usepackage{float} 
\usetikzlibrary{arrows,shapes,positioning}
\usetikzlibrary{calc,decorations.markings}
\tikzstyle{block} = [rectangle, draw, rounded corners, text centered, minimum height = 2em, minimum width = 7em, align=left, scale = 0.75]

\tikzstyle{line} = [draw, thick]

\usetikzlibrary{calc,decorations.markings}
\usepackage{hyperref}
\usepackage{amsthm}
\usepackage{makecell}
\usepackage{enumitem}
\usepackage{graphicx}

\usepackage{subfig}
\usepackage{marginnote}

\usepackage{amssymb}
\usepackage{epstopdf}
\usepackage{amsmath}
\usepackage{amsfonts}
\theoremstyle{plain}
\newtheorem{theorem}{Theorem}[section]

\theoremstyle{definition}

\theoremstyle{plain}

 \newcommand{\tor}{\text{ or }}
 \newcommand{\tand}{\text{ and }}

\def\doncl{\quad\Leftrightarrow\quad}

\newcommand{\dt}{\mathrm{d}t  }

\newcommand{\ds}{\mathrm{d}s  }

\newcommand{\HH}{\mathbb{H}}

\newcommand{\Prob}{\mathbb{P}}
\newcommand{\Proba}[1]{\mathbb{P}\left[#1\right]}

\newcommand{\branchmat}[1]{ \left\{\begin{matrix}
 #1
\end{matrix}\right. }

\def\SLEkk#1/{$\mathrm{SLE}_{#1}$}
\def\SLEk/{\SLEkk{\kappa}/}
\def\SLE/{$\mathrm{SLE}$}

\newcommand{\Z}{\mathbb{Z}}

\newcommand{\R}{\mathbb{R}}

\newcommand{\Ito}{It\^{o} }

\newcommand{\set}[1]{\left \{#1 \right\}}

\author{Tomas Kojar}

\title{A Left passage explorer model}

\begin{document}
\maketitle
{\let\thefootnote\relax\footnote{{for comments and questions please contact tkojarmathphysics@gmail.com}}}
\begin{figure}[h]
\vspace{-3.5em}
\begin{center}
\includegraphics[width=0.8\textwidth]{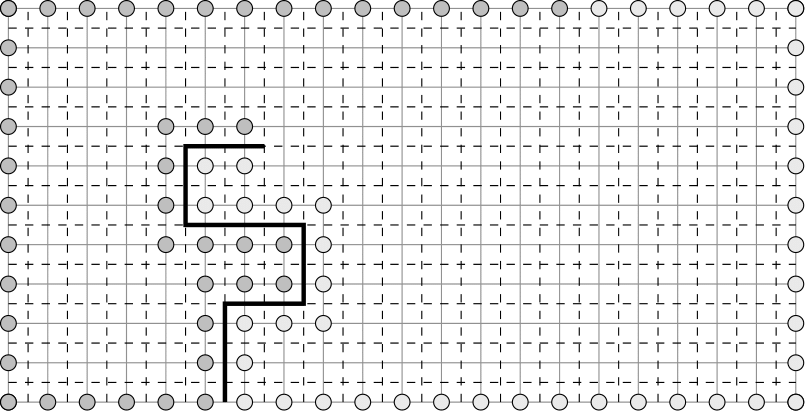}
\end{center}
\vspace{-2em}
\end{figure}

\begin{abstract} This is a short note describing a model generalizing the Harmonic explorer \cite{schramm2005harmonic} that might be of interest and it is \textit{not} intended for publication in a journal. The conjectured continuous model should have the same left-passage probability as \SLEk/ but possibly be a different curve.
\end{abstract}

\section{Introduction}
In the article \cite{schramm2005harmonic}, they constructed a model that converges to the SLE curve for $\kappa=4$. Later in their work in \cite{schramm2009contour}, they showed that the chordal contour lines of the discrete Gaussian free field converge to forms of SLE(4); after proving a height gap result, they then used techniques developed in the "toy" case of the harmonic explorer.\\
This note initially started in motivation of developing a family discrete models for each $\kappa\in (0,4)$ each of which will converge to \SLEk/. So a natural move was to extend the harmonic explorer model. There the key observable is the left passage probability of the harmonic explorer path $\gamma$, namely the function $h(z):=\Proba{\gamma\text{ passes to the left of} z }$ and this function has nice analytical properties, namely it is a harmonic function that is equal to $\pi$ on the right side of $\gamma$ and $0$ on the left.\\

\paragraph{Acknowledgements:}I thank Ilya Binder for his help in defining the model.

\section{The Left-Passage pde (LP-pde)}
 Schramm had developed in \cite{schramm2001percolation} a formula for the left passage probability for \SLEk/ for each $\kappa\in [0,8)$.
\begin{theorem}\label{whichside}
Let $\kappa\in[0,8)$, and let $z_0=x_0+i\,y_0\in\mathbb{H}$. Then the trace $\gamma$ of chordal \SLEk/ satisfies
$$
h(z_{0}):=\Proba{\gamma\hbox{ passes to the left of }z_0}=
\frac 12 +\frac{\Gamma(4/\kappa)} {\sqrt{\pi}\,\Gamma\bigl(\frac{8-\kappa}{2\kappa}\bigr)}\, \frac{x_0}{y_0}\,  F_{2,1}\Bigl(\frac 12,\frac 4\kappa,\frac 32,-\frac {x_0^2}{y_0^2}\Bigr)\,.$$,
where $F_{2,1}$ is the hypergeometric function.
 \end{theorem}
Similarly in \cite[section 5]{lawler2009conformal},  Lawler obtains a formula that is reminiscent of the argument function from the harmonic explorer:
\begin{theorem}
Let $h(z)$ again denote the probability that z is on the left side of the \SLEk/ path $\gamma(0,\infty)$, then \cite[section 5]{lawler2009conformal}
\begin{equation*}
h(z)=\int^{arg(z)}\frac{1}{sin(\theta)^{2-\frac{8}{\kappa}} }d\theta.    
\end{equation*}
\end{theorem}
% In the process of deriving this and using scale invariance he obtained the following ODE
% \begin{equation*}
%   \frac{\kappa}{2}h''(w)+\frac{42}{w^{2}+1}h'(w)=0, w\in(a,b), \twith  h(a)=0 \tand h(b)=1. 
% \end{equation*}
Using either this formula or starting with the ODEs derived in their proofs, it is easy to show that the function $f(arg(x+iy))=f(\theta)=\int^{\theta}_{0}\dfrac{1}{sin(s)^{2-\frac{8}{\kappa}}  }\ds$ will satisfy the degenerate elliptic pde:
\begin{equation}\label{LPPDE}
2(\frac{4}{\kappa}-1)\frac{\partial f}{\partial y} +y\Delta f=0\doncl div(y^{\beta}\nabla u)=0,\tag{LP-pde}
\end{equation}
where $\beta:=2(\frac{4}{\kappa}-1)$. We see that in the case $\kappa=4$ we obtain the harmonic function observable for the harmonic explorer. So starting from here it is natural to study a generalized harmonic explorer built out of the observable satisfying \eqref{LPPDE}. (One should note that this pde doesn't satisfy conformal invariance for $\kappa\neq 4$ and so this fact probably will show that this model has no relation with \SLEk/).\\
The formulation of the harmonic explorer model and its the convergence to the limit, is done studying the two-dimensional random walk $(R_{1},R_{2})$ studied in the domain $D\setminus \gamma[0,n]$ where $\gamma$ is the harmonic explorer up to time $n$. So naturally, we study the diffusion corresponding to the \eqref{LPPDE}. Because of the degeneracy of the pde we need the following Feynman-Kac theorem from \cite{feehan2015martingale}. \\
They generally consider a time-dependent, degenerate-elliptic differential operator defined by \emph{unbounded} coefficients $(a,b)$ on the half-space $\HH := \R^{d-1}\times(0,\infty)$ with $d\geq 1$,
\begin{equation}
\label{eq:MartingaleGenerator}
\frac{1}{2}\sum_{i,j=1}^d x_da_{ij}(t,x)v_{x_ix_j}(x) + \sum_{i=1}^d b_i(t,x)v_{x_i}(x), \quad (t,x) \in [0,\infty)\times\HH,
\end{equation}
and $a=(a_{ij})$, $b=(b_i)$, and $v \in C^{2}(\overline\HH)$. In our case $a_{ii}=2, a_{i,j}=0, i\neq j$ and $b_{1}=0,b_{2}=\beta$. So it satisfies the list of assumptions in their "assumptions 2.2". So we get the correspondence to the following degenerated SDE system 
\begin{equation*}
\branchmat{dX_{1,t}=2\sqrt{X_{2.t}}dW_{1,t}\\dX_{2,t}=2\sqrt{X_{2.t}}dW_{2,t}+\beta\dt}.
\end{equation*}
By \Ito and time change we get
\begin{equation*}
\branchmat{dX_{1,t}=dW_{1,t}\\dX_{2,t}=dW_{2,t}+\frac{1}{X_{2,t}}\frac{\beta-1}{2}\dt}.
\end{equation*}

So this a two dimensional tupple of (Brownian motion, Bessel process). For the Brownian motion part the natural discretization is the random walk. For the Bessel part we use the result in \cite{csaki2009transient}. In that paper they consider a nearest neighbor (NN) random walk, defined as follows: let $X_0=0,\ X_1,X_2,\ldots$ be a Markov
chain with
\begin{eqnarray}\label{defff}
E_i:=&\Proba{(X_{n+1}=i+1\mid X_n=i}=1-\Proba{X_{n+1}=i-1\mid X_n=i}\\
=&\left\{\begin{array}{ll} 1\quad & {\rm if}\quad  i=0\\\nonumber 1/2+p_i\quad & {\rm if}\quad i=1,2,\ldots,\end{array}\right.
\end{eqnarray}
where $-1/2\leq p_i\leq 1/2,\ i=1,2,\ldots$. In case $0< p_i\leq
1/2$ the sequence $\{X_i\}$ describes the motion of a particle which
starts at zero, moves over the nonnegative integers and going away
from 0 with a larger probability than to the direction of 0. They take $p_i\sim B/4i$ with $B>0$ as $i\to\infty$. Then they show that in certain sense, this Markov chain is a discrete analogue of continuous Bessel process and establish a strong invariance principle between these two processes. In particular, consider
$Y_\nu(t),\, t\geq 0$, a Bessel process of order $\nu$,
$Y_\nu(0)=0$, and let $X_n,\, n=0,1,2,\ldots$ be an NN random
walk with  $p_0=p_1=1/2$,
\begin{equation}
p_R=\frac{(R-1)^{-2\nu}-R^{-2\nu}}{(R-1)^{-2\nu}-(R+1)^{-2\nu}}-\frac12,\qquad
R=2,3,\ldots \label{pr}
\end{equation}
One main result is a strong invariance principle concerning Bessel
process and NN random walk 
\begin{theorem}
On a suitable probability space they construct a Bessel process
$\{Y_\nu(t),\, t\geq 0\},$ $\nu>0$ and an NN random walk $\{X_n,\,
n=0,1,2,\ldots\}$ with $p_R$ as in {\rm (\ref{pr})} such that for
any $\varepsilon>0$, as $n\to\infty$ we have
\begin{equation}
Y_\nu(n)-X_n=O(n^{1/4+\varepsilon})\qquad {\rm a.s.} \label{inv}
\end{equation}
\end{theorem}

\section{The LP-model}
\par
\begin{figure}[H]
    \centering
\includegraphics{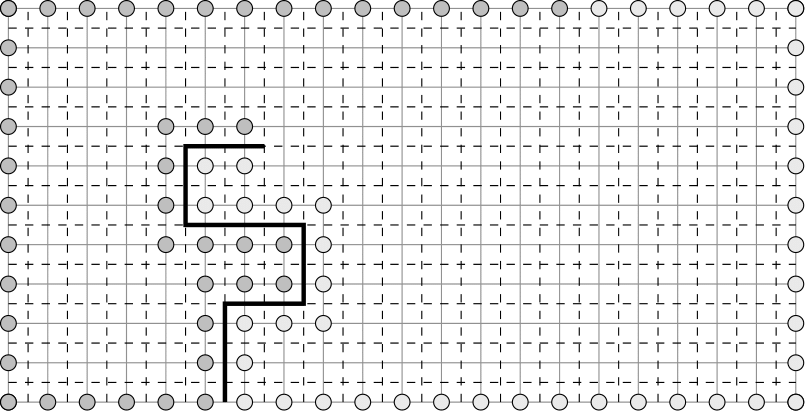}
\caption{Entire rectangle D with boundary conditions}
\label{entire_domain}
\end{figure}

 We consider lattice $\delta \Z^{2} \cap \mathbb{H}$  ,for some $\delta>0$, and rectangle D with bottom side centered at the origin.  We also include its medial lattice, which in the figure below we colored its lines with dashed lines and its vertices with black dots. Finally, for fixed boundary vertices  $v_{0},v_{end}$,  we impose boundary conditions 0 (darkgray) and 1 (lightgray) on separate segments of $\partial D=\partial_{0}D\sqcup \partial_{1}D$, see figure \ref{entire_domain} (dark-gray on left and light gray on the right respectively). 
\par The model is similar to the harmonic explorer model in \cite{schramm2005harmonic} with the exception that we replace the two dimensional random walk on the hexagonal grid by running a random walk along the x-coordinate and a discretization of the Bessel process along the y-coordinate, which we denote by (R,B). In particular, we assign to each vertex $w$ on the original lattice (black dots),  the value 
\begin{equation*}
\mathfrak{h}_{n}(w):=\Prob_{w}[ (R_{T_{D_{n}}},B_{T_{D_{n}}})\in  \partial_{1}D_{n}  ], 
\end{equation*}
where $T_{D_{n}}$ is the exit time of (R,B) from domain $D_{n}:=D\setminus \gamma[0,n]$, for path $\gamma$ as described below. The exploration path $\gamma$ runs along the medial lattice from $v_{0}$ to $v_{end}$. It moves in three posible directions: left, straight or right. We will try to use similar notation as in \cite[section 3.1]{schramm2005harmonic} as much as possible even though here the lattice is square and not triangular. \\
\begin{figure}[H]
    \centering
\includegraphics{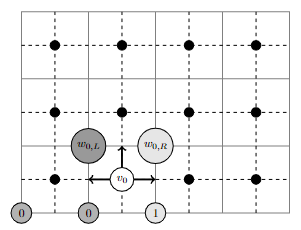}
\caption{The first step of the path}
\label{firststep}
\end{figure}
First let $V_{0}$ the set of vertices in $\partial D$, $A^{-}$ the dark-gray/negative-vertices in $V_0$ going from $v_{0}$ to $v_{end}$ but on the left and $A^{+}$ the light-gray/positive ones on the right also going from $v_{0}$ to $v_{end}$.\\
For the first-step we let $h_{0}:V_{0}\to \set{0,1}$ to be 0 on $V_{0}\cap A^{-}:=\partial_{0} D$ and to be 1 on $V_{0}\cap A^{+}:=\partial_{1} D$. On the original lattice we have two vertices $w_{0,L}$ and $w_{0,R}$ on the square containing $v_0$ (see figure \ref{firststep}). We let $p_{0,L},p_{0,R}$ equal the value of the probabilities 
\begin{equation*}
p_{0,L}:=\mathfrak{h}_{0}(w_{0,L})=\Prob_{w_{0,L}}[ (R_{T_{D}},B_{T_{D}})\in  \partial_{1}D  ]\tand p_{0,R}:=\mathfrak{h}_{0}(w_{0,R})=\Prob_{w_{0,R}}[ (R_{T_{D}},B_{T_{D}})\in  \partial_{1}D  ].
\end{equation*}
This is the analogous step of considering the harmonic extension of $h_{0}$. From here there are two possible variations.
\paragraph{variation 1}
If one could show that $p_{0,L}\leq p_{0,R}$, then as in the Harmonic explorer we consider an iid sequence of $\set{X_{n}}_{n\geq 1}$ with uniform distribution Unif([0,1]). Then we decide whether the path moves left,straight or right depending on whether:
\begin{equation*}
  X_{1}\leq p_{0,L}   ,X_{1}\in (p_{0,L},p_{0,R}) \tor p_{0,R}\leq X_{1}.
\end{equation*}
The statement $p_{0,L}\leq p_{0,R}$ is clear by a path-reflection/swapping argument. Suppose there is a trajectory $G(k):=(B_{k},R_{k})_{k=0}^{T_{D}}$  with $G(0)=w_{0,L}$ that exits on the right side $\partial_{1}D$. Then a different trajectory $\widetilde{G}(k)$ starting with $\widetilde{G}(0)=w_{0,R}$ either exists through $\partial_{1}D$ earlier or is \textit{forced} to intersect with $G(k)$ at some point $\tau$ and from then on it is identified with it till their exit at $\partial_{1}D$. Therefore, we have an inclusion of paths.

\paragraph{variation 2}
We consider two independent sequences $\set{X_{n,L}}_{n\geq 1}$,$\set{X_{n,R}}_{n\geq 1}$ to be iid with uniform distribution Unif([0,1]) to help us decide which of the three directions to pick. Then We make our choice starting clockwise from $w_{0,L}$ i.e. by by first flipping the $X_{1,L}$-coin. If $X_{1,L}\leq p_{0,L}$, we go left and if $X_{1,L}\geq p_{0,L}$ we will move either straight or right of $w_{0,R}$. So next we flip the second coin $X_{1,R}$ to help us i.e. if $X_{1,R}\leq p_{0,R}$ we move left of $w_{0,R}$ and if $X_{1,R}\geq p_{0,R}$ we move right of it.
\begin{figure}[H]
    \centering
\includegraphics{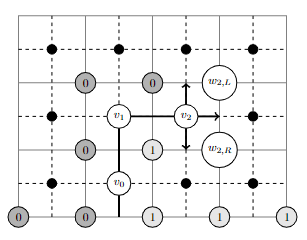}
\caption{The first two steps of the path}
\end{figure}
After the first step is done, we extend the boundary to $V_{1}:=V_{0}\cup \set{w_{0,L},w_{0,R}}$ and assign labels 0 or 1 depending on whether the path moved left, straight or right. We let $D_{1}:=D\setminus \gamma[0,1]$ denote the new domain and $\partial_{0} D_{1},\partial_{1} D_{1}$ its new labeled boundaries. For the next steps we similarly decide whether to move on the left, straight or right by using $p_{2,i}:=\mathfrak{h}_{1}(w_{2,i})=\Prob_{w_{2},i}[ (R_{T_{D_{1}}},B_{T_{D_{1}}})\in  \partial_{1}D_{1}  ] $ for $i=L,R$.

\section{Questions}
Here are some questions
\begin{enumerate}
    
    \item If shown to converge to a continuous limit, this model might serve as a counterexample of a model having having same left-passage probability as \SLEk/ but not being related to it or even having conformal invariance.

    \item Continuing on the program in \cite{schramm2009contour}, it would be interesting to consider a possible coupling of those continuous limits with the Generalized Gaussian field whose covariance is the Green function for the divergence operator $div(y^{\beta}\nabla u)$ (eg. see the fields studied in \cite{gu2017generalized}).

    \item It would be interesting to further generalize this model to other 2d-discrete processes such as for the discretization of two general \Ito diffusions.
    
\end{enumerate}

\bibliographystyle{plain} 
\bibliography{\jobname}

\begin{thebibliography}{1}

\bibitem{csaki2009transient}
Endre Cs{\'a}ki, Ant{\'o}nia F{\"o}ldes, and P{\'a}l R{\'e}v{\'e}sz.
\newblock Transient nearest neighbor random walk and bessel process.
\newblock {\em Journal of Theoretical Probability}, 22(4):992--1009, 2009.

\bibitem{feehan2015martingale}
Paul Feehan and Camelia Pop.
\newblock On the martingale problem for degenerate-parabolic partial
  differential operators with unbounded coefficients and a mimicking theorem
  for ito processes.
\newblock {\em Transactions of the American Mathematical Society},
  367(11):7565--7593, 2015.

\bibitem{gu2017generalized}
Yu~Gu and Jean-Christophe Mourrat.
\newblock On generalized gaussian free fields and stochastic homogenization.
\newblock {\em Electronic Journal of Probability}, 22:1--21, 2017.

\bibitem{lawler2009conformal}
Gregory Lawler.
\newblock Conformal invariance and 2 statistical physics.
\newblock {\em Bulletin of the American Mathematical Society}, 46(1):35--54,
  2009.

\bibitem{schramm2001percolation}
Oded Schramm.
\newblock A percolation formula.
\newblock {\em Electronic Communications in Probability}, 6:115--120, 2001.

\bibitem{schramm2005harmonic}
Oded Schramm and Scott Sheffield.
\newblock Harmonic explorer and its convergence to sle4.
\newblock {\em The Annals of Probability}, 33(6):2127--2148, 2005.

\bibitem{schramm2009contour}
Oded Schramm and Scott Sheffield.
\newblock Contour lines of the two-dimensional discrete gaussian free field.
\newblock {\em Acta mathematica}, 202(1):21--137, 2009.

\end{thebibliography}

\end{document}